\documentclass[11pt,a4paper]{article}

\usepackage{theorem,enumerate}
\usepackage{amsmath,latexsym,amssymb,amsfonts}
\usepackage{eucal}
\usepackage{color}
\usepackage{comment}

\newcounter{Scounter}
\theorembodyfont{\normalfont\slshape}
\newtheorem{thm}{Theorem}[section]

\newtheorem{Thm}{Theorem}

\newtheorem{prop}[thm]{Proposition}
\newtheorem{lem}[thm]{Lemma}

\newtheorem{Q}[thm]{Question} 

\newtheorem{claim}{Claim}[section] 

\newtheorem{fact}[claim]{Fact} 
\newtheorem{con}{Conjecture}

\numberwithin{equation}{section}
\newtheorem{remark}{\normalfont\itshape Remark}

\newcommand{\proof}{\medbreak\noindent\textit{Proof.}\quad}
\newcommand{\qed}{{$\quad\square$\vs{3.6}}}

\newcommand{\vs}[1]{\vspace*{#1 mm}}

\def\B{{ \mathcal{B}}}
\def\C{{ \mathcal{C}}}

\def\F{{ \mathcal{F}}}

\def\H{{ \mathcal{H}}}

\bfseries\normalfont

\title{Sufficient conditions for the existence of a path-factor which are related to odd components}

\author{
Yoshimi Egawa$^1$ \and\
Michitaka Furuya$^1$\footnote{\texttt{michitaka.furuya@gmail.com}} \and\
Kenta Ozeki$^2$ $^3$\footnote{\texttt{ozeki@nii.ac.jp}} \vs{5}\\
$^1$\textsl{Department of Mathematical Information Science,} \\
\textsl{Tokyo University of Science,}\\
\textsl{1-3 Kagurazaka, Shinjuku-ku, Tokyo 162-8601, Japan }\\
$^2$\textsl{National Institute of Informatics,} \\
\textsl{2-1-2 Hitotsubashi, Chiyoda-ku, Tokyo 101-8430, Japan}\\
$^3$\textsl{JST, ERATO, Kawarabayashi Large Graph Project, Japan} 
}

\date{}

\begin{document}

\maketitle

\begin{abstract}
In this paper, we are concerned with sufficient conditions for the existence of a $\{P_{2},P_{2k+1}\}$-factor.
We prove that for $k\geq 3$, there exists $\varepsilon _{k}>0$ such that if a graph $G$ satisfies $\sum _{0\leq j\leq k-1}c_{2j+1}(G-X)\leq \varepsilon _{k}|X|$ for all $X\subseteq V(G)$, then $G$ has a $\{P_{2},P_{2k+1}\}$-factor, where $c_{i}(G-X)$ is the number of components $C$ of $G-X$ with $|V(C)|=i$.
On the other hand, we construct infinitely many graphs $G$ having no $\{P_{2},P_{2k+1}\}$-factor such that $\sum _{0\leq j\leq k-1}c_{2j+1}(G-X)\leq \frac{32k+141}{72k-78}|X|$ for all $X\subseteq V(G)$.
\end{abstract}



\section{Introduction}\label{sec1}

In this paper, we consider only finite undirected simple graphs.
Let $G$ be a graph.
We let $V(G)$ and $E(G)$ denote the vertex set and the edge set of $G$, respectively.
For $u\in V(G)$, we let $N_{G}(u)$ and $d_{G}(u)$ denote the {\it neighborhood} and the {\it degree} of $u$, respectively.
For $U\subseteq V(G)$, we let $N_{G}(U)=(\bigcup _{u\in U}N_{G}(u))-U$.
For a subgraph $H$ of $G$ and a set $X\subseteq V(G)$, we let $H[X]$ denote the subgraph of $H$ induced by $V(H)\cap X$.
A graph is {\it odd} if its order is odd.
We let $\C(G)$ and $\C_{\rm odd}(G)$ denote the set of components of $G$ and the set of odd components of $G$, respectively.
Set $c(G)=|\C(G)|$ and $c_{\rm odd}(G)=|\C_{\rm odd}(G)|$.
For two graphs $H_{1}$ and $H_{2}$, we let $H_{1}\cup H_{2}$ and $H_{1}+H_{2}$ denote the {\it union} and the {\it join} of $H_{1}$ and $H_{2}$, respectively.
For a graph $H$ and an integer $s\geq 2$, we let $sH$ denote the union of $s$ disjoint copies of $H$.
Let $K_{n}$, $P_{n}$ and $C_{n}$ denote the {\it complete graph}, the {\it path} and the {\it cycle} of order $n$, respectively.
For terms and symbols not defined here, we refer the reader to \cite{D}.

Let again $G$ be a graph.
A subset $M$ of $E(G)$ is a {\it matching} if no two distinct edges in $M$ have a common endvertex.
If there is no fear of confusion, we often identify a matching $M$ of $G$ with the subgraph of $G$ induced by $M$.
A matching $M$ of $G$ is {\it perfect} if $V(M)=V(G)$.
If $G-u$ has a perfect matching for every $u\in V(G)$, $G$ is called {\it hypomatchable}.
For a set $\F$ of connected graphs, a spanning subgraph $F$ of $G$ is called an {\it $\F$-factor} if each component of $F$ is isomorphic to a graph in $\F$.
Note that a perfect matching can be regarded as a $\{P_{2}\}$-factor.
A $\{P_{n}:n\geq 2\}$-factor of $G$ is called a {\it path-factor} of $G$.
Since every path of order at least $2$ can be partitioned into paths of orders $2$ and $3$, a graph has a path-factor if and only if it has a $\{P_{2},P_{3}\}$-factor.
Akiyama, Avis and Era~\cite{AAE} gave a necessary and sufficient condition for the existence of a path-factor (here $i(G)$ denotes the number of isolated vertices of a graph $G$).

\begin{Thm}[Akiyama, Avis and Era~\cite{AAE}]
\label{ThmA}
A graph $G$ has a $\{P_{2},P_{3}\}$-factor if and only if $i(G-X)\leq 2|X|$ for all $X\subseteq V(G)$.
\end{Thm}

On the other hand, it follows from a result of Loebal and Poljak~\cite{LP} that for $k\geq 2$, the existence problem of a $\{P_{2},P_{2k+1}\}$-factor is {\bf NP}-complete.
Thus we are interested in a useful sufficient condition for the existence of a $\{P_{2},P_{2k+1}\}$-factor (for detailed historical background and motivations, we refer the reader to \cite{EF1}).

In order to state our results, we need some more preparations.
For $j\geq 1$, let $\C _{j}(H)$ be the set of components $C$ of a graph $H$ with $|V(C)|=j$, and set $c_{j}(H)=|\C _{j}(H)|$.
Note that $c_{1}(H)$ is the number of isolated vertices of $H$ (i.e., $c_{1}(H)=i(H)$).
Since no odd graph of order at most $2k-1$ has a $\{P_{2},P_{2k+1}\}$-factor, the existence of an odd subgraph of order at most $2k-1$ can be regarded as an obstacle to the existence of a $\{P_{2},P_{2k+1}\}$-factor.
Furthermore, for $k\geq 1$, if a graph $G$ has a $\{P_{2},P_{2k+1}\}$-factor, then $\sum _{0\leq j\leq k-1}(k-j)c_{2j+1}(G-X)\leq (k+1)|X|$ for all $X\subseteq V(G)$ (see \cite[Proposition~2.1]{EF2}).
Thus if a condition concerning $c_{2j+1}(G-X)~(0\leq j\leq k-1)$ for $X\subseteq V(G)$ assures us the existence of a $\{P_{2},P_{2k+1}\}$-factor, then it will make a useful sufficient condition.

Recently, Egawa and Furuya~\cite{EF1,EF2} began such a study, and they proved the following theorems.

\begin{Thm}[Egawa and Furuya~\cite{EF1}]
\label{ThmB}
Let $G$ be a graph.
If $c_{1}(G-X)+\frac{2}{3}c_{3}(G-X)\leq \frac{4}{3}|X|$ for all $X\subseteq V(G)$, then $G$ has a $\{P_{2},P_{5}\}$-factor.
\end{Thm}

\begin{Thm}[Egawa and Furuya~\cite{EF2}]
\label{ThmC}
Let $G$ be a graph.
If $c_{1}(G-X)+\frac{1}{3}c_{3}(G-X)+\frac{1}{3}c_{5}(G-X)\leq \frac{2}{3}|X|$ for all $X\subseteq V(G)$, then $G$ has a $\{P_{2},P_{7}\}$-factor.
\end{Thm}

\begin{Thm}[Egawa and Furuya~\cite{EF2}]
\label{ThmD}
Let $G$ be a graph.
If $c_{1}(G-X)+c_{3}(G-X)+\frac{2}{3}c_{5}(G-X)+\frac{1}{3}c_{7}(G-X)\leq \frac{2}{3}|X|$ for all $X\subseteq V(G)$, then $G$ has a $\{P_{2},P_{9}\}$-factor.
\end{Thm}

They also showed that the coefficients of $|X|$ in the above theorems are best possible.

These results naturally suggest the following problem:
For $k\geq 5$, is there a number $\varepsilon _{k}>0$ such that if a graph $G$ satisfies $\sum _{0\leq j\leq k-1}c_{2j+1}(G-X)\leq \varepsilon _{k}|X|$ for all $X\subseteq V(G)$, then $G$ has a $\{P_{2},P_{2k+1}\}$-factor?
Our first result in this paper is the following theorem, which gives an affirmative solution to the problem.

\begin{thm}
\label{thm1}
Let $k\geq 3$ be an integer, and let $G$ be a graph.
If $\sum _{0\leq j\leq k-1}c_{2j+1}(G-X)\leq \frac{5}{6k^{2}}|X|$ for all $X\subseteq V(G)$, then $G$ has a $\{P_{2},P_{2k+1}\}$-factor.
\end{thm}

In \cite{EF1}, Egawa and Furuya constructed examples which show that for $k\geq 3$ with $k\equiv 0~(\mbox{mod }3)$, there exist infinitely many graphs $G$ having no $\{P_{2},P_{2k+1}\}$-factor such that $\sum _{0\leq j\leq k-1}c_{2j+1}(G-X)\leq \frac{4k+6}{8k+3}|X|+\frac{2k+3}{8k+3}$ for all $X\subseteq V(G)$, and proposed the following conjecture.

\begin{con}
\label{conj1}
Let $k\geq 3$ be an integer, and let $G$ be a graph.
If $\sum _{0\leq j\leq k-1}c_{2j+1}(G-X)\leq \frac{4k+6}{8k+3}|X|$ for all $X\subseteq V(G)$, then $G$ has a $\{P_{2},P_{2k+1}\}$-factor.
\end{con}

Theorems~\ref{ThmC} and \ref{ThmD} imply that Conjecture~\ref{conj1} is true for $k\in \{3,4\}$.
Note also that Conjecture~\ref{conj1}, if true, would give an affirmative solution to the problem mentioned in the paragraph preceding Theorem~\ref{thm1} with $\varepsilon _{k}=\frac{4k+6}{8k+3}$.
However, the fact is that Conjecture~\ref{conj1} is false for large integers $k$.
Our second result is the following.

\begin{thm}
\label{thm2}
For $k\geq 29$, there exist infinitely many graphs $G$ having no $\{P_{2},P_{2k+1}\}$-factor such that $\sum _{0\leq j\leq k-1}c_{2j+1}(G-X)\leq \frac{32k+141}{72k-78}|X|$ for all $X\subseteq V(G)$.
\end{thm}

For $k\geq 36$, by simple calculations, we have $\frac{32k+141}{72k-78}<\frac{4k+6}{8k+3}$.
This implies that the coefficient of $|X|$ in Conjecture~\ref{conj1} is irrelevant for large $k$.

In Section~\ref{sec-suf}, we give a sufficient condition for the existence of an $\F$-factor for a set $\F$ with $P_{2}\in \F$.
In Section~\ref{sec-hypno}, we study fundamental properties of hypomatchable graphs without $\{P_{2},P_{2k+1}\}$-factors.
By using results in Sections~\ref{sec-suf} and \ref{sec-hypno}, we prove Theorem~\ref{thm1} in Section~\ref{sec-pf1}.
In Section~\ref{sec-coex}, we construct graphs which show that Theorem~\ref{thm2} holds.
We remark that Lemmas~\ref{lem-hypno-prep-1} and \ref{lem-hypno-prep-2}, which are proved in Section~\ref{sec-hypno}, hold for hypomatchable graphs in general, and thus could hopefully be used in the study of other types of factors.

In our proof, we make use of the following facts.

\begin{fact}
\label{fact1}
Let $k\geq 2$ be an integer, and let $G$ be a graph.
Then $G$ has a $\{P_{2},P_{2k+1}\}$-factor if and only if $G$ has a path-factor $F$ such that $\C_{2i+1}(F)=\emptyset $ for every $i~(1\leq i\leq k-1)$.
\end{fact}

\begin{fact}
\label{fact2}
Let $a,b,c,d\in \mathbb{R}$ with $c\not=0$ and $ad-bc\geq 0$.
Then the function $f(x)=\frac{ax+b}{cx+d}$ is non-decreasing in the interval $x>-\frac{d}{c}$.
\end{fact}

\section{A sufficient condition for the existence of a component-factor}\label{sec-suf}

Let $\F$ be a set of connected graphs.
For a graph $H$, we let $\B_{\F}(H)$ denote the set of those hypomatchable components of $H$ which have no $\F$-factor, and set $b_{\F}(H)=|\B_{\F}(H)|$.
Cornu\'{e}jols and Hartvigsen~\cite{CH} proved that when $P_{2}\in \F$ and $\F-\{P_{2}\}$ consists of hypomatchable graphs, a graph $G$ has an $\F$-factor if and only if $b_{\F}(G-X)\leq |X|$ for all $X\subseteq V(G)$.
The following proposition follows from the proof of the ``if'' part of the above result of Cornu\'{e}jols and Hartvigsen, but we include its proof for the convenience of the reader.

\begin{prop}
\label{prop-suf-1}
Let $\F$ be a set of connected graphs such that $P_{2}\in \F$, and let $G$ be a graph.
If $b_{\F}(G-X)\leq |X|$ for all $X\subseteq V(G)$, then $G$ has an $\F$-factor.
\end{prop}

In our proof of Proposition~\ref{prop-suf-1}, we choose a set $S$ of vertices of $G$ so that
\begin{enumerate}
\item[{\bf (S1)}]
$c_{\rm odd}(G-S)-|S|$ is as large as possible, and
\item[{\bf (S2)}]
subject to (S1), $|S|$ is as large as possible.
\end{enumerate}
Note that $c_{\rm odd}(G-S)-|S|\geq c_{\rm odd}(G)-|\emptyset |\geq 0$ (it is possible that $S=\emptyset $, but our argument in this section works even if $S=\emptyset $).

We make use of the following lemma, which was proved in \cite{EF2}.

\begin{lem}[Egawa and Furuya~\cite{EF2}]
\label{lem-suf-1}
Let $G$ be a graph, and let $S$ be a subset of $V(G)$ satisfying (S1) and (S2).
Then the following hold.
\begin{enumerate}[{\upshape(i)}]
\item
We have $\C(G-S)=\C_{\rm odd}(G-S)$.
\item
For each $C\in \C_{\rm odd}(G-S)$, $C$ is hypomatchable.
\item
Let $H$ be the bipartite graph with bipartition $(S,\C_{\rm odd}(G-S))$ defined by letting $uC\in E(H)~(u\in S,C\in \C_{\rm odd}(G-S))$ if and only if $N_{G}(u)\cap V(C)\not=\emptyset $.
Then for every $X\subseteq S$, $|N_{H}(X)|\geq |X|$.
\end{enumerate}
\end{lem}

\medbreak\noindent\textit{Proof of Proposition~\ref{prop-suf-1}.}\quad
Let $G$ be as in Proposition~\ref{prop-suf-1}.
Choose $S\subseteq V(G)$ so that (S1) and (S2) hold.
Set $T=\C(G-S)$, and let $H$ be the bipartite graph $H$ with bipartition $(S,T)$ defined by letting $uC\in E(H)~(u\in S,C\in T)$ if and only if $N_{G}(u)\cap V(C)\not=\emptyset $.
By Lemma~\ref{lem-suf-1}(i)(ii), each element of $T$ is hypomatchable.
Let $T_{1}=\{C\in T:C$ has no $\F$-factor$\}$.

\begin{claim}
\label{cl-thm-suf-1}
For every $Y\subseteq T_{1}$, $|N_{H}(Y)|\geq |Y|$.
\end{claim}
\proof
Suppose that there exists $Y\subseteq T_{1}$ such that $|N_{H}(Y)|<|Y|$.
Set $X'=N_{H}(Y)$.
Then each element of $Y$ is a hypomatchable component of $G-X'$ having no $\F$-factor, and hence $|Y|\leq b_{\F}(G-X')$.
Consequently $|X'|=|N_{H}(Y)|<|Y|\leq b_{\F}(G-X')$, which contradicts the assumption of the theorem.
\qed

It follows from Claim~\ref{cl-thm-suf-1} and Hall's marriage theorem that $H$ has a matching covering $T_{1}$.
Let $M$ be a maximum matching of $H$ covering $T_{1}$.

We show that $M$ covers $S$ by using an alternating-path argument.
Suppose that $S-V(M)\not=\emptyset $.
Let $v\in S-V(M)$.
An {\it alternating path} is a path of $H$ starting from $v$ and alternately containing edges in $E(H)-M$ and edges in $M$.
Let $A$ be the set of those vertices of $H$ which are contained in an alternating path.
By the definition of an alternating path, every vertex in $A\cap S$ except $v$ belongs to $V(M)$ and, for $x\in V(M)\cap S$, $x$ belongs to $A$ if and only if the other endvertex of the edge in $M$ which is incident with $x$ belongs to $A\cap T$.
Thus $|A\cap S|=|A\cap V(M)\cap S|+1=|A\cap V(M)\cap T|+1$.
By the definition of an alternating path, we also have $N_{H}(A\cap S)\subseteq A\cap T$.
Therefore it follows from Lemma~\ref{lem-suf-1}(iii) that $A\cap T\not\subseteq A\cap V(M)\cap T$.
Take $u\in (A\cap T)-(A\cap V(M)\cap T)$, and let $P$ be an alternating path connecting $v$ and $u$.
Then $M'=(M-E(P))\cup (E(P)-M)$ is a matching of $H$ which covers $T_{1}$ and satisfies $|V(M')|=|V(M)\cup \{v,u\}|>|V(M)|$, which contradicts the maximality of $M$.
Consequently $M$ covers $S\cup T_{1}$.

Recall that each element of $T$ is a hypomatchable graph.
Thus for $uC\in M~(u\in S,C\in T)$, the subgraph of $G$ induced by $\{u\}\cup V(C)$ has a perfect matching.
Since each element of $V(H)-V(M)~(\subseteq T-T_{1})$ has an $\F$-factor, it follows that $G$ has an $\F$-factor.

This completes the proof of Proposition~\ref{prop-suf-1}.
\qed

\section{Hypomatchable graphs having no $\{P_{2},P_{2k+1}\}$-factor}\label{sec-hypno}

For an integer $k\geq 1$ and a set $\F$ of connected graphs with $P_{2}\in \F$, a pair $(\varepsilon ,\lambda )~(\varepsilon >0,\lambda \in \mathbb{N})$ is {\it $(k,\F)$-good} if the following holds:
every hypomatchable graph $G$ of order at least $2k+1$ with no $\F$-factor has a set $X\subseteq V(G)$ with $|X|\geq \lambda $ such that $\sum _{0\leq j\leq k-1}c_{2j+1}(G-X)\geq \varepsilon |X|$.


In this section, we study the existence of a $(k,\{P_{2},P_{2k+1}\})$-good pair.
In Subsection~\ref{sec-hypno-ear}, we state fundamental properties of odd ear decompositions of hypomatchable graphs.
In Subsection~\ref{sec-hypno-prep}, we introduce several notions related to odd ear decompositions, and prove two lemmas which we use in Subsection~\ref{sec-hypno-pf}.
In Subsection~\ref{sec-hypno-pf}, we show that there exists a $(k,\{P_{2},P_{2k+1}\})$-good pair for each $k\geq 3$ by proving the following proposition.

\begin{prop}
\label{prop-hypno-1}
Let $k\geq 3$ be an integer.
Then $\left(\frac{1}{k^{2}},5\right)$ is a $(k,\{P_{2},P_{2k+1}\})$-good pair.
\end{prop}

\subsection{Odd ear decompositions for hypomatchable graphs}\label{sec-hypno-ear}

We start with a structure theorem for hypomatchable graphs.
Let $G$ be a graph.
A sequence $(H_{1},\ldots ,H_{m})$ of edge-disjoint subgraphs of $G$ is an {\it odd ear decomposition} if
\begin{enumerate}
\item[{\bf (E1)}]
$V(G)=\bigcup _{1\leq i\leq m}V(H_{i})$;
\item[{\bf (E2)}]
for each $1\leq i\leq m$, $|E(H_{i})|$ is odd and $|E(H_{i})|\geq 3$;
\item[{\bf (E3)}]
$H_{1}$ is a cycle; and
\item[{\bf (E4)}]
for each $2\leq i\leq m$, either
\begin{enumerate}
\item[{\bf (E4-1)}]
$H_{i}$ is a path and only the endvertices of $H_{i}$ belong to $\bigcup _{1\leq j\leq i-1}V(H_{j})$, or
\item[{\bf (E4-2)}]
$H_{i}$ is a cycle with $|V(H_{i})\cap (\bigcup _{1\leq j\leq i-1}V(H_{j}))|=1$.
\end{enumerate}
\end{enumerate}

Lov\'{a}sz~\cite{L} proved the following theorem.

\begin{Thm}[Lov\'{a}sz~\cite{L}]
\label{Thm-hypno-ear-A}
Let $G$ be a graph with $|V(G)|\geq 3$.
Then $G$ is hypomatchable if and only if $G$ has an odd ear decomposition.
\end{Thm}

By observing the proof of Theorem~\ref{Thm-hypno-ear-A}, we obtain the following theorem.

\begin{Thm}[Lov\'{a}sz~\cite{L}]
\label{Thm-hypno-ear-B}
Let $G$ be a hypomatchable graph, and let $G_{0}$ be a subgraph of $G$.
If $G_{0}$ has an odd ear decomposition $\H=(H_{1},\ldots ,H_{m})$ and $G-V(G_{0})$ has a perfect matching, then $\H$ can be extended to an odd ear decomposition $(H_{1},\ldots ,H_{m},H_{m+1},\ldots ,H_{m'})$ of $G$.
\end{Thm}

In \cite{EF2}, the following lemma was proved.

\begin{lem}[Egawa and Furuya~\cite{EF2}]
\label{lem-hypno-ear-1}
Let $G$ be a hypomatchable graph, and let $(H_{1},\ldots ,H_{m})$ be an odd ear decomposition of $G$.
Then for each $i~(2\leq i\leq m)$, there exists an odd ear decomposition $(H'_{1},\ldots ,H'_{m'})$ of $G$ such that $H_{i}\subseteq H'_{1}$.
\end{lem}

\subsection{Height and related definitions}\label{sec-hypno-prep}

Let $G$ be a hypomatchable graph of order at least three, and let $\H=(H_{1},\ldots ,H_{m})$ be an odd ear decomposition of $G$.
We assume that we have chosen $\H$ so that
\begin{enumerate}
\item[{\bf (H1)}]
$(|V(H_{1})|,\ldots ,|V(H_{m})|)$ is lexicographically as large as possible.
\end{enumerate}

For each $i~(1\leq i\leq m)$, let $Q(i)=H_{i}-\bigcup _{1\leq j\leq i-1}V(H_{j})$.
Note that $V(Q(i))\cap V(H_{j})=\emptyset $ for any $i,j$ with $i>j$, and $\bigcup _{1\leq j\leq i}V(H_{j})=\bigcup _{1\leq j\leq i}V(Q(j))$ for each $i$.
We have $Q(1)=H_{1}$ and, by (E2) and (E4), $Q(i)$ is a path of even order for $2\leq i\leq m$.

\begin{lem}
\label{lem-hypno-prep-1}
Let $G$ be a hypomatchable graph of order at least three.
Let $(H_{1},\ldots ,H_{m})$ be an odd ear decomposition of $G$ satisfying (H1), and write $|V(H_{1})|=2l+1$.
Then for each $i~(2\leq i\leq m)$, we have $|V(Q(i))|\leq 2l$, where $Q(i)$ is as defined above.
\end{lem}
\proof
Suppose that $|V(Q(i))|\geq 2l+1$.
Since $|V(Q(i))|$ is even, this forces $|V(Q(i))|\geq 2l+2$.
By Lemma~\ref{lem-hypno-ear-1}, there exists an odd ear decomposition $(H'_{1},\ldots ,H'_{m'})$ of $G$ such that $H_{i}\subseteq H'_{1}$.
Then $|V(H'_{1})|\geq |V(Q(i))|\geq 2l+2>|V(H_{1})|$, which contradicts (H1).
Thus $|V(Q(i))|\leq 2l$.
\qed

Now for $1\leq i\leq m$ and $x\in V(Q(i))$, we recursively define the {\it height} ${\rm ht}(x)$ of $x$, the {\it height} ${\rm ht}(H_{i})$ of $H_{i}$, the set $I(x)$ of indices, and the path $R(x)$ as follows.

For each $x\in V(Q(1))$, let ${\rm ht}(x)=0$ and $I(x)=\{1\}$, and let $R(x)$ be a spanning path of $H_{1}$ with an endvertex $x$.
Let ${\rm ht}(H_{1})=0$.

Let $2\leq i\leq m$, and assume that we have defined ${\rm ht}(y)$, ${\rm ht}(H_{j})$, $I(y)$ and $R(y)$ for all $1\leq j\leq i-1$ and $y\in V(Q(j))$.
Take $x\in V(Q(i))$.
Then there exist two edge-disjoint paths $Q$ and $Q'$ on $H_{i}$ connecting $x$ and $\bigcup _{1\leq j\leq i-1}V(H_{j})$.
Since $E(H_{i})$ is odd, precisely one of $Q$ and $Q'$ has even length (i.e., odd order).
Let $H_{i}(x)$ denote the one which has odd order, and $y_{x}$ denote the endvertex of $H_{i}(x)$ different from $x$.
Note that $y_{x}\in \bigcup _{1\leq j\leq i-1}V(H_{j})$.
Define ${\rm ht}(x)={\rm ht}(y_{x})+1$ and $I(x)=I(y_{x})\cup \{i\}$.
Let $R(x)$ be the path defined by $R(x)=H_{i}(x)\cup R(y_{x})$.
Let ${\rm ht}(H_{i})=\min\{{\rm ht}(y):y\in V(Q(i))\}$.

\begin{claim}
\label{cl-hypno-prep-1}
For $i~(1\leq i\leq m)$ and $x\in V(Q(i))$, the following hold:
\begin{enumerate}
\item[{\upshape(i)}]
$\{1,i\}\subseteq I(x)\subseteq \{1,\ldots ,i\}$;
\item[{\upshape(ii)}]
${\rm ht}(x)=|I(x)|-1$;
\item[{\upshape(iii)}]
for $j\in I(x)-\{i\}$, ${\rm ht}(x)>{\rm ht}(H_{j})$;
\item[{\upshape(iv)}]
for $j~(1\leq j\leq m)$, $j\in I(x)$ if and only if $V(R(x))\cap V(Q(j))\not=\emptyset $;
\item[{\upshape(v)}]
$V(H_{1})\subseteq V(R(x))$; and
\item[{\upshape(vi)}]
for $j\in I(x)-\{1\}$, both $|V(R(x))\cap V(Q(j))|$ and $|V(Q(j))-V(R(x))|$ are even and $|V(R(x))\cap V(Q(j))|\geq 2$.
\end{enumerate}
\end{claim}
\proof
We proceed by induction on $i$.
If $i=1$, then the claim clearly holds.
Thus let $2\leq i\leq m$, and assume that all $1\leq j\leq i-1$ and $y\in V(Q(j))$ satisfy (i)--(vi).

Let $j_{0}$ be the index such that $y_{x}\in V(Q(j_{0}))$.
By the induction assumption, $j_{0}$ and $y_{x}$ satisfy (i)--(vi).
Since $I(x)=I(y_{x})\cup \{i\}$ and $\{1\}\subseteq I(y_{x})\subseteq \{1,\ldots ,j_{0}\}$, (i) holds.
Since ${\rm ht}(y_{x})=|I(y_{x})|-1$, we have
\begin{align*}
{\rm ht}(x)={\rm ht}(y_{x})+1=|I(y_{x})|=|I(x)-\{i\}|=|I(x)|-1,
\end{align*}
which implies (ii).
Since ${\rm ht}(y_{x})>{\rm ht}(H_{j})$ for $j\in I(y_{x})-\{j_{0}\}$ and ${\rm ht}(y_{x})\geq \min\{{\rm ht}(y):y\in V(Q(j_{0}))\}={\rm ht}(H_{j_{0}})$,
\begin{align*}
{\rm ht}(x)>{\rm ht}(y_{x})\geq {\rm ht}(H_{j})\mbox{~~~for all }j\in I(y_{x})~(=I(x)-\{i\}),
\end{align*}
and hence (iii) holds.
Since $j\in I(y_{x})~(=I(x)-\{i\})$ if and only if $V(R(y_{x}))\cap V(Q(j))\not=\emptyset $, it follows from $R(x)=H_{i}(x)\cup R(y_{x})$ that (iv) holds.
We have $V(H_{1})\subseteq V(R(y_{x}))\subseteq V(R(x))$, which implies (v).
Now we show (vi).
For $j\in I(x)-\{1,i\}$, $|V(R(y_{x}))\cap V(Q(j))|~(=|V(R(x))\cap V(Q(j))|)$ and $|V(Q(j))-V(R(y_{x}))|~(=|V(Q(j))-V(R(x))|)$ are even and $|V(R(x))\cap V(Q(j))|\geq 2$.
Since $V(R(x))\cap V(Q(i))=V(H_{i}(x))-\{y_{x}\}$ and $|V(Q(i))|$ is even, both $|V(R(x))\cap V(Q(i))|$ and $|V(Q(i))-V(R(x))|$ are even.
Since $x\in V(R(x))\cap V(Q(i))$, this implies that $|V(R(x))\cap V(Q(i))|\geq 2$.
Thus (vi) holds.
\qed

For $i~(2\leq i\leq m)$, since $|V(Q(i))|$ is even, $\{{\rm ht}(x):x\in V(Q(i))\}=\{{\rm ht}(x):x$ is an endvertex of $V(Q(i))\}$, and hence $Q(i)$ has an endvertex $u_{i}$ such that ${\rm ht}(u_{i})={\rm ht}(H_{i})$.

The following lemma is the main result of this subsection.

\begin{lem}
\label{lem-hypno-prep-2}
Let $G$ be a hypomatchable graph of order at least three, and let $(H_{1},\ldots ,H_{m})$ be an odd ear decomposition of $G$ satisfying (H1).
Let $h_{0}$ be an integer with $1\leq h_{0}\leq \max\{{\rm ht}(H_{i}):1\leq i\leq m\}$.
Then the set $\{u_{i}:2\leq i\leq m,{\rm ht}(H_{i})=h_{0}\}$ is an independent set of  $G$ where, as above, $u_{i}$ denotes an endvertex of $Q(i)$ such that ${\rm ht}(u_{i})={\rm ht}(H_{i})$.
\end{lem}
\proof
Suppose that $\{u_{i}:2\leq i\leq m,{\rm ht}(H_{i})=h_{0}\}$ is not an independent set of $G$.
Then there exist two indices $i$ and $i'$ with $2\leq i<i'\leq m$ such that ${\rm ht}(H_{i})={\rm ht}(H_{i'})=h_{0}$ and $u_{i}u_{i'}\in E(G)$.
By Claim~\ref{cl-hypno-prep-1}(ii), $h_{0}={\rm ht}(H_{i'})={\rm ht}(u_{i'})=|I(u_{i'})|-1$.
If $i\in I(u_{i'})$, then ${\rm ht}(H_{i'})={\rm ht}(u_{i'})>{\rm ht}(H_{i})$ by Claim~\ref{cl-hypno-prep-1}(iii), which contradicts the fact that ${\rm ht}(H_{i})={\rm ht}(H_{i'})$.
Thus
\begin{align}
i\not\in I(u_{i'}).\label{eq-cl-hypno-pf-6-1}
\end{align}

Let $R_{1}$ be the subpath on $H_{i}$ with $V(R_{1})\supseteq V(Q(i))$ connecting $u_{i}$ and $y_{u_{i}}$.
Note that $1\in I(u_{i'})$.
Let $R_{2}$ be the shortest subpath on $R(u_{i'})$ connecting $u_{i'}$ and $\bigcup _{1\leq j\leq i-1}V(H_{j})$.
Let $H^{*}=(R_{1}\cup R_{2})+u_{i}u_{i'}$.
By (\ref{eq-cl-hypno-pf-6-1}), $H^{*}$ is a path or a cycle.
Furthermore, if $H^{*}$ is a path, then only the endvertices of $H^{*}$ belong to $\bigcup _{1\leq j\leq i-1}V(H_{j})$; if $H^{*}$ is a cycle, then $|V(H^{*})\cap (\bigcup _{1\leq j\leq i-1}V(H_{j}))|=|\{y_{u_{i}}\}|=1$.
Since both $|E(R_{1})|$ and $|E(R_{2})|$ are even by Claim~\ref{cl-hypno-prep-1}(vi), $|E(H^{*})|$ is odd.
In particular, $(H_{1},\ldots ,H_{i-1},H^{*})$ is an odd ear decomposition of the subgraph of $G$ induced by $(\bigcup _{1\leq j\leq i-1}V(H_{j}))\cup V(H^{*})$.
By Claim~\ref{cl-hypno-prep-1}(vi), $|V(Q(j))-V(H^{*})|=|V(Q(j))-V(R_{2})|$ is even for every $j\in I(u_{i'})$ with $j\geq i+1$.
Note that $V(Q(i))-V(H^{*})=V(Q(i))-V(R_{1})=\emptyset $.
Since $|V(Q(j)|$ is even for every $j$ with $2\leq j\leq m$, it now follows from Claim~\ref{cl-hypno-prep-1}(iv) that $|V(Q(j))-V(H^{*})|$ is even for every $j$ with $i\leq j\leq m$.
Consequently $G-((\bigcup _{1\leq j\leq i-1}V(H_{j}))\cup V(H^{*}))$ has a perfect matching.
Therefore by Theorem~\ref{Thm-hypno-ear-B}, $G$ has an odd ear decomposition $(H_{1},\ldots ,H_{i-1},H^{*},H'_{1},\ldots ,H'_{m'})$.
Since $|V(H^{*})|=(|V(H_{i})|-1)+|V(R_{2})|\geq |V(H_{i})|-1+|V(Q(i'))|+1>|V(H_{i})|$, this contradicts the assumption (H1).
\qed

\subsection{Proof of Proposition~\ref{prop-hypno-1}}\label{sec-hypno-pf}

Let $k\geq 3$, and let $G$ be a hypomatchable graph of order at least $2k+1$ having no $\{P_{2},P_{2k+1}\}$-factor.
We use the notation introduced in the preceding subsection.
In particular, we choose an odd ear decomposition $(H_{1},\ldots ,H_{m})$ of $G$ so that (H1) holds and, for each $1\leq i\leq m$, let $u_{i}$ denote an endvertex of $Q(i)$ such that ${\rm ht}(u_{i})={\rm ht}(H_{i})$.
Having Lemma~\ref{lem-hypno-prep-2} in mind, we aim at showing that there exists an integer $h_{1}$ with $1\leq h_{1}\leq \max\{{\rm ht}(H_{i}):1\leq i\leq m\}$ for which $|\{i:1\leq i\leq m,{\rm ht}(H_{i})=h_{1}\}|$ is ``large'' (see Claim~\ref{cl-hypno-pf-4}).

A set $I\subseteq \{1,2,\ldots ,m\}$ of indices with $1\in I$ is {\it admissible} if the subgraph of $G$ induced by $\bigcup _{i\in I}V(Q(i))$ has a $\{P_{2},P_{2k+1}\}$-factor.

\begin{claim}
\label{cl-hypno-pf-1}
There is no admissible set.
\end{claim}
\proof
Suppose that there exists an admissible set $I$.
Then the subgraph of $G$ induced by $\bigcup _{i\in I}V(Q(i))$ has a $\{P_{2},P_{2k+1}\}$-factor $F$.
On the other hand, for each $i$ with $2\leq i\leq m$ and $i\not\in I$, from the fact that $Q(i)$ is a path of even order, we see that $Q(i)$ has a perfect matching $M_{i}$.
Since $\{V(Q(i)):i\not\in I\}$ is a partition of $V(G)-(\bigcup _{i\in I}V(Q(i)))$, $F\cup (\bigcup _{i\not\in I}M_{i})$ is a $\{P_{2},P_{2k+1}\}$-factor of $G$, which is a contradiction.
\qed

Write $|V(H_{1})|=2l+1$ and $|V(G)|=2n+1$.
Set $h=\max\{{\rm ht}(H_{i}):1\leq i\leq m\}$.
Let $a_{j}=\min\{k-l-j,l\}$ for each integer $j\geq 1$, and let $a^{*}=\sum _{1\leq j\leq h}a_{j}$.
We now prove three claims.

\begin{claim}
\label{cl-hypno-pf-2}
For $i~(2\leq i\leq m)$, $|V(Q(i))|\leq 2a_{{\rm ht}(H_{i})}$.
\end{claim}
\proof
In view of Lemma~\ref{lem-hypno-prep-1}, it suffices to show that $|V(Q(i))|\leq 2(k-l-{\rm ht}(H_{i}))$.
By Claim~\ref{cl-hypno-prep-1}(i)(ii), $|I(u_{i})-\{1,i\}|=|I(u_{i})|-2={\rm ht}(u_{i})-1$.
By Claim~\ref{cl-hypno-prep-1}(vi), $|V(R(u_{i}))\cap V(Q(j))|\geq 2$ for each $j\in I(u_{i})-\{1,i\}$.
Furthermore, by Claim~\ref{cl-hypno-prep-1}(v) and the definition of $u_{i}$ and $R(u_{i})$, we have $V(H_{1})\cup V(Q(i))\subseteq V(R(u_{i}))$.
Hence $|V(R(u_{i}))|\geq |V(H_{1})|+2|I(u_{i})-\{1,i\}|+|V(Q(i))|=(2l+1)+2({\rm ht}(u_{i})-1)+|V(Q(i))|=2l+2{\rm ht}(u_{i})-1+|V(Q(i))|$.
Note that $|V(R(u_{i}))|$ is odd by Claim~\ref{cl-hypno-prep-1}(iv)(v)(vi).
If $|V(R(u_{i}))|\geq 2k+1$, then by Claim~\ref{cl-hypno-prep-1}(vi) and Fact~\ref{fact1}, $I(u_{i})$ is an admissible set, which contradicts Claim~\ref{cl-hypno-pf-1}.
Thus $|V(R(u_{i}))|\leq 2k-1$.
Consequently $2k-1\geq |V(R(u_{i}))|\geq 2l+2{\rm ht}(u_{i})-1+|V(Q(i))|$.
This implies that $|V(Q(i))|\leq 2k-2l-2{\rm ht}(u_{i})=2(k-l-{\rm ht}(H_{i}))$, as desired.
\qed

\begin{claim}
\label{cl-hypno-pf-3}
We have $h\leq k-l-1$.
\end{claim}
\proof
Let $i_{0}~(1\leq i_{0}\leq m)$ be an index such that ${\rm ht}(H_{i_{0}})=h$.
Then by Claim~\ref{cl-hypno-pf-2}, $2\leq |V(Q(i_{0}))|\leq 2(k-l-{\rm ht}(H_{i_{0}}))$.
Thus $h={\rm ht}(H_{i_{0}})\leq k-l-1$.
\qed

Note that since $\{1\}$ and $\{1,2\}$ are not admissible by Claim~\ref{cl-hypno-pf-1}, we have
\begin{align}
|V(H_{1})|+|V(Q(2))|\leq 2k-1.\label{eq-hypno-pf-0}
\end{align}
This implies $m\geq 3$, and hence $h\geq 1$.

\begin{claim}
\label{cl-hypno-pf-4}
There exists $h_{1}~(1\leq h_{1}\leq h)$ such that $|\{i:1\leq i\leq m,{\rm ht}(H_{i})=h_{1}\}|\geq \frac{n-l}{a^{*}}$.
\end{claim}
\proof
For each $j~(1\leq j\leq h)$, let $N_{j}=\{i:1\leq i\leq m, {\rm ht}(H_{i})=j\}$.
Let $h_{1}~(1\leq h_{1}\leq h)$ be an integer such that $|N_{h_{1}}|=\max\{|N_{j}|:1\leq j\leq h\}$.
We show that $h_{1}$ is a desired integer.
It follows from Claim~\ref{cl-hypno-pf-2} that
\begin{align*}
(2n+1)-(2l+1) &= |V(G)|-|V(H_{1})|\\
&=\sum _{2\leq i\leq m}|V(Q(i))|\\
&\leq \sum _{2\leq i\leq m}2a_{{\rm ht}(H_{i})}\\
&= \sum _{1\leq j\leq h}\left(\sum _{i\in N_{j}}2a_{{\rm ht}(H_{i})}\right)\\
&= \sum _{1\leq j\leq h}\left(\sum _{i\in N_{j}}2a_{j}\right)\\
&= \sum _{1\leq j\leq h}2a_{j}|N_{j}|\\
&\leq \sum _{1\leq j\leq h}2a_{j}|N_{h_{1}}|\\
&=2a^{*}|N_{h_{1}}|.
\end{align*}
Consequently $|\{i:1\leq i\leq m,{\rm ht}(H_{i})=h_{1}\}|=|N_{h_{1}}|\geq \frac{n-l}{a^{*}}$.
\qed

We can now complete the proof of Proposition~\ref{prop-hypno-1}.
Let $h_{1}$ be as in Claim~\ref{cl-hypno-pf-4}, and set $X=V(G)-\{u_{i}:2\leq i\leq m,{\rm ht}(H_{i})=h_{1}\}$.
Then it follows from Lemma~\ref{lem-hypno-prep-2} and Claim~\ref{cl-hypno-pf-4} that $\sum _{0\leq j\leq k-1}c_{2j+1}(G-X)=c_{1}(G-X)=|\{u_{i}:2\leq i\leq m,{\rm ht}(H_{i})=h_{1}\}|\geq \frac{n-l}{a^{*}}$ and $|X|=2n+1-|\{u_{i}:2\leq i\leq m,{\rm ht}(H_{i})=h_{1}\}|\leq 2n+1-\frac{n-l}{a^{*}}$.
Since
\begin{align*}
\frac{n-l}{a^{*}}-\frac{n-l}{2na^{*}+a^{*}-n+l}|X| &\geq \frac{n-l}{a^{*}}-\frac{n-l}{2na^{*}+a^{*}-n+l}\left(2n+1-\frac{n-l}{a^{*}}\right)=0,
\end{align*}
we obtain
\begin{align}
\sum _{0\leq j\leq k-1}c_{2j+1}(G-X)\geq \frac{n-l}{a^{*}}\geq \frac{n-l}{2na^{*}+a^{*}-n+l}|X|.\label{eq-hypno-pf-1}
\end{align}
Since $m\geq 3$, we also have $|X|\geq |V(H_{1})|+(m-1)\geq 5$.

We give a rough bound for $n$, $l$ and $a^{*}$.
Since $2n+1=|V(G)|\geq 2k+1$, we have $n\geq k$.
By (\ref{eq-hypno-pf-0}), $(2l+1)+2\leq |V(H_{1})|+|V(Q(2))|\leq 2k-1$, and hence $l\leq k-2$.
We show that $a^{*}\leq \frac{(k-2)(k-1)}{2}$.
By Claim~\ref{cl-hypno-pf-3} and the definition of $a_{i}$ and $a^{*}$, $a^{*}=\sum _{1\leq j\leq h}a_{j}\leq \sum _{1\leq j\leq h}(k-l-j)\leq \sum _{1\leq j\leq k-l-1}(k-l-j)$.
Since $l\geq 1$, it follows that $a^{*}\leq \sum _{1\leq j\leq k-2}(k-1-j)=\frac{(k-2)(k-1)}{2}$.

By (\ref{eq-hypno-pf-1}),
\begin{align*}
\sum _{0\leq j\leq k-1}c_{2j+1}(G-X)&\geq \frac{n-l}{2na^{*}+a^{*}-n+l}|X|\\
&\geq \frac{n-(k-2)}{2n\frac{(k-2)(k-1)}{2}+\frac{(k-2)(k-1)}{2}-n+(k-2)}|X|\\
&= \frac{2n-2k+4}{(2k^{2}-6k+2)n+(k^{2}-k-2)}|X|.
\end{align*}
Recall that $k\geq 3$.
By Fact~\ref{fact2}, the function $\frac{2x-2k+4}{(2k^{2}-6k+2)x+(k^{2}-k-2)}~(x\geq k)$ is non-decreasing.
If $n>\frac{1}{2}k^{2}$, then $\frac{2n-2k+4}{(2k^{2}-6k+2)n+(k^{2}-k-2)}\geq \frac{2\cdot \frac{1}{2}k^{2}-2k+4}{(2k^{2}-6k+2)\frac{1}{2}k^{2}+(k^{2}-k-2)}>\frac{1}{k^{2}}$; if $n\leq \frac{1}{2}k^{2}$, then since $\sum _{0\leq j\leq k-1}c_{2j+1}(G-X)\geq 1$ and $|X|\leq |V(G)|-1=2n\leq k^{2}$, we clearly have $\sum _{0\leq j\leq k-1}c_{2j+1}(G-X)\geq \frac{1}{k^{2}}|X|$.
Thus $\sum _{0\leq j\leq k-1}c_{2j+1}(G-X)\geq \frac{1}{k^{2}}|X|$.

Therefore the set $X$ satisfies $|X|\geq 5$ and $\sum _{0\leq j\leq k-1}c_{2j+1}(G-X)\geq \frac{1}{k^{2}}|X|$.
Since $G$ is arbitrary, this completes the proof of Proposition~\ref{prop-hypno-1}.

\begin{remark}
\label{re-hypno-pf-1}
{\rm
We used rough estimates for $n$, $l$ and $a^{*}$ in the proof of Proposition~\ref{prop-hypno-1} because our aim was to show the existence of a $k$-good pair.
If we go through some more calculations, we will get a $k$-good pair $(\varepsilon ,\lambda )$ with a larger value of $\varepsilon $ than $\frac{1}{2k^{2}}$.
}
\end{remark}

\section{Proof of Theorem~\ref{thm1}}\label{sec-pf1}

In view of Proposition~\ref{prop-hypno-1}, Theorem~\ref{thm1} immediately follows from the following proposition.

\begin{prop}
\label{prop-pf1-1}
Let $k\geq 1$ be an integer and $\F$ be a set of connected graphs with $P_{2}\in \F$, and let $(\varepsilon ,\lambda )$ be a $(k,\F)$-good pair with $\varepsilon \leq 1$.
If a graph $G$ satisfies $\sum _{0\leq j\leq k-1}c_{2j+1}(G-X)\leq \frac{\lambda \varepsilon }{\lambda +1}|X|$ for all $X\subseteq V(G)$, then $G$ has an $\F$-factor.
\end{prop}
\proof
Suppose that $G$ has no $\F$-factor.
Then by Proposition~\ref{prop-suf-1}, there exists $X'\subseteq V(G)$ such that $b_{\F}(G-X')>|X'|$.
Let $\B_{1}=\{C\in \B_{\F}(G-X'):|V(C)|\geq 2k+1\}$.
Then by the definition of a $(k,\F)$-good pair, each $C\in \B_{1}$ has a set $X_{C}\subseteq V(C)$ with $|X_{C}|\geq \lambda $ such that $\sum _{0\leq j\leq k-1}c_{2j+1}(C-X_{C})\geq \varepsilon |X_{C}|$.
Set $X_{0}=X\cup (\bigcup _{C\in \B_{1}}X_{C})$.
Then
$$
\C_{2j+1}(G-X_{0})=\C_{2j+1}(G-X')\cup \left(\bigcup _{C\in \B_{1}}\C_{2j+1}(C-X_{C})\right)\mbox{ for }0\leq j\leq k-1.
$$
Consequently
\begin{align*}
\sum _{0\leq j\leq k-1}c_{2j+1}(G-X_{0})&=\sum _{0\leq j\leq k-1}c_{2j+1}(G-X')+\sum _{0\leq j\leq k-1}\left(\sum _{C\in \B_{1}}c_{2j+1}(C-X_{C})\right)\\
&=\sum _{0\leq j\leq k-1}c_{2j+1}(G-X')+\sum _{C\in \B_{1}}\left(\sum _{0\leq j\leq k-1}c_{2j+1}(C-X_{C})\right)\\
&\geq \sum _{0\leq j\leq k-1}c_{2j+1}(G-X')+\varepsilon \sum _{C\in \B_{1}}|X_{C}|.
\end{align*}
This together with the assumption that $\varepsilon \leq 1$ and the fact that $|X_{C}|\geq \lambda ~(C\in \B_{1})$ implies
\begin{align*}
&\sum _{0\leq j\leq k-1}c_{2j+1}(G-X_{0})-\frac{\lambda \varepsilon }{\lambda +1}\left(\sum _{0\leq j\leq k-1}c_{2j+1}(G-X')+\sum _{C\in \B_{1}}(|X_{C}|+1)\right)\\
&\geq \sum _{0\leq j\leq k-1}c_{2j+1}(G-X')+\varepsilon \sum _{C\in \B_{1}}|X_{C}|-\frac{\lambda \varepsilon }{\lambda +1}\left(\sum _{0\leq j\leq k-1}c_{2j+1}(G-X')+\sum _{C\in \B_{1}}(|X_{C}|+1)\right)\\
&= \left(1-\frac{\lambda \varepsilon }{\lambda +1}\right)\sum _{0\leq j\leq k-1}c_{2j+1}(G-X')+\frac{\varepsilon }{\lambda +1}\sum _{C\in \B_{1}}(|X_{C}|-\lambda )\\
&\geq 0,
\end{align*}
and hence
$$
\frac{\lambda \varepsilon }{\lambda +1}\left(\sum _{0\leq j\leq k-1}c_{2j+1}(G-X')+\sum _{C\in \B_{1}}(|X_{C}|+1)\right)\leq \sum _{0\leq j\leq k-1}c_{2j+1}(G-X_{0}).
$$
Therefore
\begin{align*}
\frac{\lambda \varepsilon }{\lambda +1}|X_{0}| &= \frac{\lambda \varepsilon }{\lambda +1}\left(|X'|+\sum _{C\in \B_{1}}|X_{C}|\right)\\
&< \frac{\lambda \varepsilon }{\lambda +1}\left(b_{\F}(G-X')+\sum _{C\in \B_{1}}|X_{C}|\right)\\
&= \frac{\lambda \varepsilon }{\lambda +1}\left((|\B_{\F}(G-X')-\B_{1}|)+\sum _{C\in \B_{1}}(|X_{C}|+1)\right)\\
&\leq  \frac{\lambda \varepsilon }{\lambda +1}\left(\sum _{0\leq j\leq k-1}c_{2j+1}(G-X')+\sum _{C\in \B_{1}}(|X_{C}|+1)\right)\\
&\leq \sum _{0\leq j\leq k-1}c_{2j+1}(G-X_{0}),
\end{align*}
which contradicts the assumption of the proposition.
\qed

\section{Proof of Theorem~\ref{thm2}}\label{sec-coex}

Throughout this section, we fix an integer $k\geq 29$.
Set $l=2\left\lfloor \frac{k-12}{17}\right\rfloor $, $m=\left\lfloor \frac{2k-l+1}{8}\right\rfloor$ and $r=2k+1-l-8m~\left(=2k+1-l-8\left\lfloor \frac{2k-l+1}{8}\right\rfloor \right)$.
Then the following lemma holds.

\begin{lem}
\label{lem-coex-1}
\begin{enumerate}
\item[{\upshape(i)}]
$2k+1=l+8m+r$.
\item[{\upshape(ii)}]
$\frac{2k-56}{17}\leq l\leq \frac{2k-24}{17}$.
\item[{\upshape(iii)}]
$\frac{16k-39}{68}\leq m\leq \frac{32k+73}{136}$.
\item[{\upshape(iv)}]
$m\geq 2l+3$.
\item[{\upshape(v)}]
$r\in \{1,3,5,7\}$.
\end{enumerate}
\end{lem}
\proof
\begin{enumerate}
\item[{\upshape(i)}]
This follows from the definition of $r$.
\item[{\upshape(ii)}]
We have $\frac{2k-56}{17}=2\left(\frac{k-12}{17}-\frac{16}{17}\right)\leq l\leq 2\cdot \frac{k-12}{17}=\frac{2k-24}{17}$.
\item[{\upshape(iii)}]
It follows from (ii) that $\frac{16k-39}{68}=\frac{1}{8}\left(2k-\frac{2k-24}{17}+1\right)-\frac{7}{8}\leq \frac{2k-l+1}{8}-\frac{7}{8}\leq m\leq \frac{2k-l+1}{8}\leq \frac{32k+73}{136}$.
\item[{\upshape(iv)}]
It follows from (ii) and (iii) that $m\geq \frac{16k-39}{68}>\frac{2(2k-24)}{17}+2\geq 2l+2$, and hence $m\geq 2l+3$.
\item[{\upshape(v)}]
Since $0=2k+1-l-8\cdot \frac{2k-l+1}{8}\leq r<2k+1-l-8\left(\frac{2k-l+1}{8}-1\right)=8$, we have $0\leq r\leq 7$.
Since $l~(=2\left\lfloor \frac{k-12}{17}\right\rfloor )$ is even, $r$ is odd.
Thus $r\in \{1,3,5,7\}$.
\qed
\end{enumerate}

\begin{figure}
\begin{center}
{\unitlength 0.1in
\begin{picture}( 11.0000,  9.0000)(  1.6000,-10.5000)
%
\special{sh 1.000}%
\special{ia 800 200 50 50  0.0000000  6.2831853}%
\special{pn 8}%
\special{ar 800 200 50 50  0.0000000  6.2831853}%
%
\special{sh 1.000}%
\special{ia 600 400 50 50  0.0000000  6.2831853}%
\special{pn 8}%
\special{ar 600 400 50 50  0.0000000  6.2831853}%
%
\special{sh 1.000}%
\special{ia 1000 400 50 50  0.0000000  6.2831853}%
\special{pn 8}%
\special{ar 1000 400 50 50  0.0000000  6.2831853}%
%
\special{sh 1.000}%
\special{ia 400 600 50 50  0.0000000  6.2831853}%
\special{pn 8}%
\special{ar 400 600 50 50  0.0000000  6.2831853}%
%
\special{sh 1.000}%
\special{ia 600 800 50 50  0.0000000  6.2831853}%
\special{pn 8}%
\special{ar 600 800 50 50  0.0000000  6.2831853}%
%
\special{sh 1.000}%
\special{ia 1000 800 50 50  0.0000000  6.2831853}%
\special{pn 8}%
\special{ar 1000 800 50 50  0.0000000  6.2831853}%
%
\special{sh 1.000}%
\special{ia 1200 600 50 50  0.0000000  6.2831853}%
\special{pn 8}%
\special{ar 1200 600 50 50  0.0000000  6.2831853}%
%
\special{sh 1.000}%
\special{ia 800 1000 50 50  0.0000000  6.2831853}%
\special{pn 8}%
\special{ar 800 1000 50 50  0.0000000  6.2831853}%
%
\special{pn 8}%
\special{pa 800 1000}%
\special{pa 400 600}%
\special{fp}%
\special{pa 400 600}%
\special{pa 800 200}%
\special{fp}%
\special{pa 800 200}%
\special{pa 1200 600}%
\special{fp}%
\special{pa 1200 600}%
\special{pa 800 1000}%
\special{fp}%
%
\special{pn 8}%
\special{pa 1000 800}%
\special{pa 600 800}%
\special{fp}%
\special{pa 600 800}%
\special{pa 600 400}%
\special{fp}%
\special{pa 600 400}%
\special{pa 1000 400}%
\special{fp}%
\special{pa 1000 400}%
\special{pa 1000 800}%
\special{fp}%
\put(2.5000,-6.0000){\makebox(0,0){$u$}}%
\put(13.5000,-6.0000){\makebox(0,0){$v$}}%
\end{picture}}%

\caption{The graph $H$}
\label{f1}
\end{center}
\end{figure}
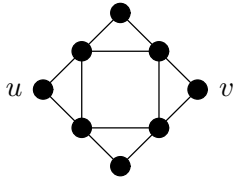

Here we construct a graph $Q$ by using an idea by Bauer, Broersma and Veldman~\cite{BBV} as follows.
Let $H$ be the graph depicted in Figure~\ref{f1} having specified vertices $u$ and $v$.
Take $m$ disjoint copies $H_{1},\ldots ,H_{m}$ of $H$ and, for each $i~(1\leq i\leq m)$, let $u_{i}$ and $v_{i}$ be the vertices of $H_{i}$ corresponding to the vertices $u$ and $v$ of $H$, respectively.
Set $U=\{u_{i},v_{i}:1\leq i\leq m\}$.
Let $R$ be a set of $r$ vertices with $(\bigcup _{1\leq i\leq m}V(H_{i}))\cap R=\emptyset $.
Let $T$ be the graph obtained from $\bigcup _{1\leq i\leq m}H_{i}$ by adding the vertices in $R$ and joining all possible pairs of vertices in $U\cup R$.
Let $L$ be a complete graph of order $l$, and let $Q=L+T$.
By Lemma~\ref{lem-coex-1}(i), $|V(Q)|=l+8m+r=2k+1$.

\begin{lem}
\label{lem-coex-3}
Let $1\leq i\leq m$, and let $X\subseteq V(H_{i})$ be a set with $\{u_{i},v_{i}\}\subseteq X$.
Then the following hold.
\begin{enumerate}
\item[{\upshape(i)}]
We have $\sum _{0\leq j\leq k-1}c_{2j+1}(H_{i}-X)\leq 2$.
\item[{\upshape(ii)}]
If $\sum _{0\leq j\leq k-1}c_{2j+1}(H_{i}-X)=1$, then $|X|\geq 3$.
\item[{\upshape(iii)}]
If $\sum _{0\leq j\leq k-1}c_{2j+1}(H_{i}-X)=2$, then $|X|\geq 4$.
\end{enumerate}
\end{lem}
\proof
Since the independence number of $H_{i}-\{u_{i},v_{i}\}$ is $2$, (i) holds.
Since $|V(H_{i})-\{u_{i},v_{i}\}|$ is even, (ii) holds.
Since $H_{i}-\{u_{i},v_{i}\}$ is $2$-connected, (iii) holds.
\qed

In view of Lemma~\ref{lem-coex-1}(iv), the following lemma follows from Theorems~3 and 4 of \cite{BBV}.

\begin{lem}[Bauer et al.~\cite{BBV}]
\label{lem-coex-4}
\begin{enumerate}
\item[{\upshape(i)}]
The graph $Q-R$ has no Hamiltonian path.
\item[{\upshape(ii)}]
For $X\subseteq V(Q)-R$, if $c(Q-(X\cup R))\geq 2$, then $c(Q-(X\cup R))\leq \frac{2m+1}{4m+l}|X|$.
\end{enumerate}
\end{lem}

\begin{lem}
\label{lem-coex-5}
The graph $Q$ has no Hamiltonian path.
\end{lem}
\proof
Suppose that $Q$ has a Hamiltonian path $P$.
By the definition of $Q$, if a vertex in $V(Q)-R$ is adjacent to a vertex $R$ on $P$, then the vertex belongs to $V(L)\cup U$.
Since $V(L)\cup U$ is a clique of $Q$, there exists a Hamiltonian path of $Q-R$ obtained from $P-R$ by adding some edges, which contradicts Lemma~\ref{lem-coex-4}(i).
\qed

\begin{lem}
\label{lem-coex-6}
For all $X\subseteq V(Q)$, $\sum _{0\leq j\leq k-1}c_{2j+1}(Q-X)\leq \frac{32k+141}{72k-78}|X|+\frac{32k+879}{288k-312}$.
\end{lem}
\proof
Since $k\geq 29$, we have $l=2\left\lfloor \frac{k-12}{17} \right\rfloor \geq 2$.
By Lemma~\ref{lem-coex-1}(iv), we also have $m\geq 7$.
Let $X\subseteq V(Q)$.
We first show that $\sum _{0\leq j\leq k-1}c_{2j+1}(Q-X)\leq \frac{2m+1}{4m+l}|X|+\frac{m-7}{8m+2l}$.

We may assume that $\sum _{0\leq j\leq k-1}c_{2j+1}(Q-X)\geq 1$.
Since $|V(Q)|=2k+1$, we have $|X|\geq 2$.
If $Q-X$ is connected, then it follows from Lemma~\ref{lem-coex-1}(iv) that $\sum _{0\leq j\leq k-1}c_{2j+1}(Q-X)=1=\frac{2m+1}{4m+l}\cdot 2+\frac{2l-4}{8m+2l}\leq \frac{2m+1}{4m+l}|X|+\frac{m-7}{8m+2l}$.
Thus we may assume that $Q-X$ is disconnected.
In particular, $V(L)\subseteq X$.

Assume first that $U\subseteq X$.
By Lemma~\ref{lem-coex-3}(i), $\sum _{0\leq j\leq k-1}c_{2j+1}(H_{i}-X)\leq 2$ for every $1\leq i\leq m$.
For each $h\in \{1,2\}$, let $m_{h}$ denote the number of $H_{i}$ such that $\sum _{0\leq j\leq k-1}c_{2j+1}(H_{i}-X)=h$.
Note that $Q-X$ have at most one component intersecting with $R$.
Hence
$$
\sum _{0\leq j\leq k-1}c_{2j+1}(Q-X)\leq \sum _{1\leq i\leq m}\left(\sum _{0\leq j\leq k-1}c_{2j+1}(H_{i}-X)\right)+1=m_{1}+2m_{2}+1.
$$
By Lemma~\ref{lem-coex-3}(ii)(iii),
$$
|X|=\sum _{1\leq i\leq m}|(V(H_{i})-\{u_{i},v_{j}\})\cap X|+|U|+|V(L)|+|R\cap X|\geq m_{1}+2m_{2}+2m+l.
$$
Therefore
\begin{align}
\sum _{0\leq j\leq k-1}c_{2j+1}(Q-X)\leq \frac{m_{1}+2m_{2}+1}{m_{1}+2m_{2}+2m+l}|X|.\label{eq-lem-coex-6.1}
\end{align}
By Fact~\ref{fact2}, $\frac{x+1}{x+2m+l}~(x\geq 0)$ is non-decreasing.
Since $m_{1}+2m_{2}\leq 2m$, it follows from (\ref{eq-lem-coex-6.1}) that
$$
\sum _{0\leq j\leq k-1}c_{2j+1}(Q-X)\leq \frac{2m+1}{4m+l}|X|<\frac{2m+1}{4m+l}|X|+\frac{m-7}{8m+2l}.
$$
Thus we may assume that $U\not\subseteq X$.
Then by the construction of $Q$, $c(Q-X)=c(Q-(X\cup R))$.
Since $Q-X$ is disconnected, $c(Q-(X\cup R))=c(Q-X)\geq 2$.
It follows from Lemma~\ref{lem-coex-4}(ii) that
\begin{align*}
\sum _{0\leq j\leq k-1}c_{2j+1}(Q-X) &\leq c(Q-X)\\
&= c(Q-(X\cup R))\\
&\leq \frac{2m+1}{4m+l}|X-R|\\
&< \frac{2m+1}{4m+l}|X|+\frac{m-7}{8m+2l}.
\end{align*}
Consequently $\sum _{0\leq j\leq k-1}c_{2j+1}(Q-X)\leq \frac{2m+1}{4m+l}|X|+\frac{m-7}{8m+2l}$.

Recall that $l\geq 2$.
By Fact~\ref{fact2}, $\frac{2x+1}{4x+l}$ and $\frac{x-7}{8x+2l}~(x>0)$ are non-decreasing.
Hence it follows from Lemma~\ref{lem-coex-1}(ii)(iii) that
$$
\frac{2m+1}{4m+l}\leq \frac{2\cdot \frac{32k+73}{136}+1}{4\cdot \frac{32k+73}{136}+l}\leq \frac{2\cdot \frac{32k+73}{136}+1}{4\cdot \frac{32k+73}{136}+\frac{2k-56}{17}}=\frac{32k+141}{72k-78}
$$
and
$$
\frac{m-7}{8m+2l}\leq \frac{\frac{32k+73}{136}-7}{8\cdot \frac{32k+73}{136}+2l}\leq \frac{\frac{32k+73}{136}-7}{8\cdot \frac{32k+73}{136}+2\cdot \frac{2k-56}{17}}=\frac{32k-879}{288k-312}.
$$
Therefore $\sum _{0\leq j\leq k-1}c_{2j+1}(Q-X)\leq \frac{32k+141}{72k-78}|X|+\frac{32k-879}{288k-312}$.
\qed

Now we are ready to prove Theorem~\ref{thm2}.
Let $n\geq 1$ be an integer.
Let $Q_{0}$ be a complete graph of order $n$.
Let $Q_{1},Q_{2},\ldots ,Q_{2n+1}$ be disjoint copies of the graph $Q$.
Let $G_{n}=Q_{0}+(\bigcup _{1\leq i\leq 2n+1}Q_{i})$.

For $1\leq i\leq 2n+1$, since $|V(Q_{i})|=2k+1$ and $Q_{i}$ has no Hamiltonian path by Lemma~\ref{lem-coex-5}, $Q_{i}$ has no $\{P_{2},P_{2k+1}\}$-factor.
Suppose that $G_{n}$ has a $\{P_{2},P_{2k+1}\}$-factor $F$.
Then for each $i~(1\leq i\leq 2n+1)$, $F$ contains an edge joining $V(Q_{i})$ and $V(Q_{0})$.
Since $2n+1>2|V(Q_{0})|$, this implies that there exists $a\in V(Q_{0})$ such that $d_{F}(a)\geq 3$, which is a contradiction.
Thus
\begin{align}
\mbox{$G_{n}$ has no $\{P_{2},P_{2k+1}\}$-factor.}\label{no-factor}
\end{align}

We next show that $\sum _{0\leq j\leq k-1}c_{2j+1}(G_{n}-X)\leq \frac{32k+141}{72k-78}|X|$ for all $X\subseteq V(G_{n})$.
Let $X\subseteq V(G_{n})$.
Assume for the moment that $V(Q_{0})\not\subseteq X$.
Then $G_{n}-X$ is connected.
Clearly we may assume that $\sum _{0\leq j\leq k-1}c_{2j+1}(G_{n}-X)=1$.
Then $|X|\geq 3$ because $|V(G_{n})|>2k+1$.
Hence $\sum _{0\leq j\leq k-1}c_{2j+1}(G_{n}-X)=1<\frac{32k+141}{72k-78}\cdot 3\leq \frac{32k+141}{72k-78}|X|$.
Thus we may assume that $V(Q_{0})\subseteq X$.
Then
\begin{align}
\C_{2j+1}(G_{n}-X)=\bigcup _{1\leq i\leq 2n+1}\C_{2j+1}(Q_{i}-X)\mbox{~~~for every }0\leq j\leq k-1.\label{eq-coex-1}
\end{align}
By Lemma~\ref{lem-coex-6} and (\ref{eq-coex-1}),
\begin{align*}
\sum _{0\leq j\leq k-1}c_{2j+1}(G_{n}-X) &= \sum _{0\leq j\leq k-1}\left(\sum _{1\leq i\leq 2n+1}c_{2j+1}(Q_{i}-X)\right)\\
&\leq \sum _{1\leq i\leq 2n+1}\left(\frac{32k+141}{72k-78}|V(Q_{i})\cap X|+\frac{32k-879}{288k-312}\right)\\
&= \frac{32k+141}{72k-78}(|X|-n)+\frac{32k-879}{288k-312}(2n+1)\\
&= \frac{32k+141}{72k-78}|X|+\frac{(-64n+32)k-2322n-879}{288k-312}\\
&< \frac{32k+141}{72k-78}|X|.
\end{align*}
Consequently
\begin{align}
\sum _{0\leq j\leq k-1}c_{2j+1}(G_{n}-X)\leq \frac{32k+141}{72k-78}|X| \mbox{ for all } X\subseteq V(G_{n}).\label{eq-coex-2}
\end{align}

By (\ref{no-factor}) and (\ref{eq-coex-2}), we obtain Theorem~\ref{thm2}.

\end{document}